\documentclass{article} 
\usepackage{graphicx} 

\newcommand{\bbibitem}{\bibitem}

\newcommand{\llabel}[1]{{\label{#1}}}

\renewcommand{\r}[1]{(\ref{#1})}
\font\tenmsb=msbm10
\font\sevenmsb=msbm7
\font\fivemsb=msbm5
\newfam\msbfam
\textfont\msbfam=\tenmsb
\scriptfont\msbfam=\sevenmsb
\scriptscriptfont\msbfam=\fivemsb
\def\Bbb#1{{\fam\msbfam\relax#1}}
\newcommand{\bi}{\begin{itemize}}
\newcommand{\ei}{\end{itemize}}
\newcommand{\bd}{\begin{description}}
\newcommand{\ed}{\end{description}}
\newcommand{\be}{\begin{enumerate}}
\newcommand{\ee}{\end{enumerate}}

\renewcommand{\i}{\item}

\newcommand{\bqn}{\begin{eqnarray}}
\newcommand{\eqn}{\end{eqnarray}}
\newcommand{\eqnn}{\nonumber\end{eqnarray}}
\newcommand{\eqnl}[1]{\llabel{#1}\end{eqnarray}}

\newcommand{\nn}{\nonumber}

\newcommand{\ba}[1]{\begin{array}{#1}}
\newcommand{\ea}{\end{array}}

\newcommand{\R}{\Bbb{R}}
\newcommand{\C}{\Bbb{C}}

\newcommand{\fine}{\end{document}}

\def \trait (#1) (#2) (#3){\vrule width #1pt height #2pt depth #3pt}
\def \qed{\hfill
        \trait (0.1) (6) (0)
        \trait (6) (0.1) (0)
        \kern-6pt   
        \trait (6) (6) (-5.9)
        \trait (0.1) (6) (0)
\medskip}
\def \qedmio{\hfill
             \trait (8) (8) (-0.1)
             \medskip}
\def \quadp{{\Huge $\qedmio$}}

\newtheorem{ml}{\bf Lemma}
\newtheorem{Theorem}{\bf Theorem}
\newtheorem{mo}{\bf \underline{{\sl Observation}}}
\newtheorem{mrem}{\bf \underline{{\sl Remark}}}
\newtheorem{mcc}{\bf Corollary}
\newtheorem{Definition}{\bf Definition}
\newtheorem{mpr}{\bf Proposition}
\newtheorem{mproperty}{\bf Property}

\newcommand{\bt}{\begin{Theorem}}
\newcommand{\et}{\end{Theorem}}
\newcommand{\bl}{\begin{ml}}
\newcommand{\el}{\end{ml}}
\newcommand{\bo}{\noindent\begin{mo}\rm}
\newcommand{\eo}{\end{mo}}
\newcommand{\bp}{\begin{mpr}}
\newcommand{\ep}{\end{mpr}}
\newcommand{\bc}{\begin{mcc}}

\newcommand{\ec}{\end{mcc}}
\newcommand{\bdeff}{\begin{Definition}}
\newcommand{\edeff}{\end{Definition}}
\newcommand{\bproperty}{\begin{mproperty}}
\newcommand{\eproperty}{\end{mproperty}}
\newcommand{\brem}{\begin{mrem}\rm}
\newcommand{\erem}{\end{mrem}}

\newcommand{\proof}{{\bf Proof. }}
\newcommand{\ppot}[4]
{

\begin{center}
~\includegraphics[height=#3truecm,width=#4truecm]{#1.eps}\\ 
{#2}
\end{center}\noindent$\!\!$}

\newcommand{\lam}{\lambda}

\newcommand{\al}{\alpha}

\newcommand{\sceq}{Schr\"{o}dinger equation\ }

\renewcommand{\k}{{\mbox{\bf k}}}
\newcommand{\p}{{\mbox{\bf p}}}

\renewcommand{\l}{\mbox{{\footnotesize \bf L}}}

\newcommand{\Sdc}{  {\cal S}_\C^d   }
\newcommand{\Tdc}{  {\cal T}_\C^d   }
\newcommand{\Suc}{  {\cal S}_\C^u   }
\newcommand{\Tuc}{  {\cal T}_\C^u   }

\newcommand{\Sdr}{  {\cal S}_\R^d   }
\newcommand{\Tdr}{  {\cal T}_\R^d   }
\newcommand{\Sur}{  {\cal S}_\R^u   }
\newcommand{\Tur}{  {\cal T}_\R^u   }

\begin{document} 
\begin{center} \noindent
{\LARGE{\sl{\bf On the $K+P$ Problem for a Three-level Quantum System:
Optimality Implies Resonance}}}
\end{center}

\vskip 1cm
\begin{center}
Ugo Boscain, Thomas Chambrion\\
{\footnotesize SISSA-ISAS, via Beirut 2-4, 34014 Trieste (Italy)}

\vspace{.5 cm}
Jean-Paul Gauthier\\
 {\footnotesize 
D\'epartement de
Math\'ematiques, Analyse Appliqu\'ee et Optimisation, 
Univerisit\'e de Bourgogne, 
9, Avenue Alain Savary B.P. 47870-21078 Dijon Cedex, France} 
\end{center}



\vspace{2.5cm} \noindent \rm 

\begin{quotation}
\noindent  {\bf Abstract}        
We apply techniques of subriemannian geometry on Lie groups to
laser-induced population transfer in a three-level quantum system. The aim
is to induce transitions
by two laser pulses, of arbitrary shape and frequency, minimizing the
pulse energy. We prove
that the
Hamiltonian system given by the 
Pontryagin Maximum Principle is completely
integrable, since this problem can be stated as a ``$\k\oplus\p$ problem''
on
a simple Lie group. Optimal trajectories and controls are exhausted.
The main
result is that optimal controls correspond to lasers that are ``in
resonance''.
\end{quotation}

\vspace{2cm}
\begin{center}
Preprint SISSA 30/2002/AF
\end{center}
\section{Introduction}
\subsection{Physical context}
In the recent past years, people started to approach the  control of the 
Schr\"odinger equation, 
using techniques of geometric control theory (see for instance 
\cite{daless,rabitz,brokko,rama}).
In this paper we apply techniques of subriemannian geometry on Lie
groups to the population transfer problem in a three-level quantum
system driven by two external fields (in the rotating wave approximation)  
of arbitrary shape and frequency.
The aim is to induce complete population transfer by minimizing the pulse
energy. 
The dynamics is governed by the time dependent Schr\"odinger equation (in
a system of units such that $\hbar=1$): 
\bqn
i\frac{d\psi(t)}{dt}=H\psi(t), 
\eqnl{se} 
where $\psi(.):\R\to\C^3$ and:
\bqn 
H=\left(\ba{ccc} E_1&\Omega_1&0\\ \Omega_1^\ast &E_2&\Omega_2\\
0&\Omega_2^\ast&E_3 \ea\right).
\eqnl{ham-gen-3}  
Here $(^\ast)$ indicates the complex conjugation involution.
The controls $\Omega_1(.),\Omega_2(.)$, that we assume  to be different 
from
zero only in a fixed interval
$[0,T]$, are 
connected to the physical parameters by
$\Omega_j(t)=\mu_{j}{{\cal F}}_j(t)/2$, $j=1,2$, with 
${\cal F}_j$ the external pulsed field and  
$\mu_{j}$ the couplings 
(intrinsic
to the quantum system) that we have restricted to
couple only levels $j$ and $j+1$ by pairs. 
\brem
This finite-dimensional problem can be thought as the reduction of an
infinite-dimensional problem in the following way. 
We start with a Hamiltonian which is the sum of a ``drift-term'' $H_0$,
plus a time dependent potential $V(t)$ (the control term, i.e. the
lasers). The drift term is assumed to be diagonal, with eigenvalues
(energy levels) 
$...>E_3>E_2>E_1.$ 
Then in this spectral resolution of
$H_0$, we assume the control term $V(t)$  to couple only the
energy levels $E_1,E_2$ and $E_2,E_3$. The projected problem in the
eigenspaces corresponding 
to $E_1,E_2,E_3$ is completely decoupled and is described by the 
Hamiltonian \r{ham-gen-3}. 
\erem
The problem is the following:\\\\
{\bf Problem.} {\it Assume that for time $t\leq0$ the state of the system 
lies in the
eigenspace corresponding to the ground eigenvalue
$E_1$. We want to determine suitable 
controls $\Omega_i(.),$ $i=1,2$, such that for time $t\geq T$,
the system reaches the
eigenspace corresponding to  $E_3$, 
requiring that these controls minimize
the cost ({\rm energy} in the following):                
\bqn 
J=\int_{0}^{T}\left(|\Omega_1(t)|^2+|\Omega_2(t)|^2\right)~dt.
\eqnl{fluence-om}}
In \cite{BCG} this problem was studied assuming that the controls
$\Omega_j$ are ``in
resonance'':
\bqn 
&&\Omega_j(t)=u_j(t)~e^{i (\omega_j t+\al_j)},~~\omega_j=E_{j+1}-E_j,\nn\\
&&u_j(.):\R\to\R,~~\al_j\in[-\pi,\pi],~~j=1,2.
\eqnl{reson}
In the sequel we call this second problem (of minimizing the cost
\r{fluence-om}, which in this case reduces to 
$\int_{0}^{T}\big((u_1(t)^2+$ $u_2(t)^2\big)~dt$),  
the {\sl
"real-resonant"} problem. The first problem (with arbitrary complex 
controls) will
be called  the {\sl "general-complex"} problem.\\
In \cite{BCG}, the  real-resonant problem was treated as follows:
\bi
\i first, using a time dependent change of coordinates that leaves
invariant the source (eigenstate 1) and the target (eigenstate 3), it was
possible to eliminate the drift term and hence to reduce the problem to a
singular-Riemannian problem over the real sphere $S^2$;
\i second, the Hamiltonian system obtained from the Pontryagin Maximum
Principle (PMP in the following)   was Liouville integrable and explicit 
expressions for geodesics
were found   (although not so simple expressions).  
\ei
In this paper, we address both problems in a more abstract setting:
\bi
\i in both cases, first, we eliminate the drift. For the 
general-complex case, we use a time dependent change of coordinates plus a
gauge
transformation. Again, this change of coordinates leaves invariant both
the
source and the target; 
\i second, in both cases, we lift the problem into a right-invariant
subriemannian control
problem on a \underline{real} simple Lie group $G$, ($G=G^{\R}=SO(3)$ for the
real-resonant problem and $G=G^{\C}=SU(3)$
for the general-complex one). 
This subriemannian problem has very special features:\\
There is an element $g_0$  of $G^{\R}\subset G^{\C}$, and subgroups 
$K^{\R}, K^{\C}$, ($K^{\R} \approx O(2), K^{\C} \approx U(2)$) 
such that, denoting by ${\k}^{\R}, {\k}^{\C}$ the Lie algebras 
of $K^{\R}, K^{\C}$ respectively:\\
\\ 
a. The Lie algebras $\l, \k$ of the pair $(G,K)$ have associated Cartan 
decomposition $\l=\k \oplus \p$, with the usual commutation relations,\\
b. The right-invariant distribution is determined by the $\p$ subspace 
of $\l$,\\
c. The right invariant subriemannian metric is determined by a scalar 
multiple of $Kil\big|_p$, where $Kil\big|_p$ is the Killing form 
restricted to $\p$,\\
d. The source  $S$ and the target $T$ of the optimal control problem 
are respectively  $S=K_{g_{0}}=  g_{0} K {g_{0}}^{-1}$, the 
conjugation of $K$ by $g_{0}$, and $T=g_{0}S$;
 
\i subriemannian problems on semi-simple Lie groups in which the
distribution is determined by the $\p$ subspace of a Cartan decomposition 
and the metric
is proportional to $Kil\big|_p$ (called in the following  
{\it $(\k\oplus\p)$-problems}) have two important features:
\bi 
\i abnormal extremals (if any) are never
optimal since the so called {\it Goh condition} (which is a necessary
condition for optimality of abnormal extremals, see \cite{AS}) is never
satisfied (see Appendix C); \i normal extremals can be computed with 
 a very powerful technique developped
by Jurdjevic in \cite{jurd-MCT,jurd-SCL}. 
In particular he proved that the
Hamiltonian system given by the PMP is completely integrable and he gave
explicit expressions
for  geodesics. For sake of completeness, the Jurdjevic's technique is 
summarized 
in Appendix B;
\ei
\i as the PMP requires, we apply  suitable transversality
conditions, corresponding to the source and the target. These conditions
allow to restrict the set of admissible extremals.
\i we prove that
all the geodesics satisfying the transversality
conditions have the same length. Optimality follows.
\ei
\brem
In the paper \cite{jurd-SCL}, the connection between the geodesics for
the $\k\oplus\p$ problem and
the geodesics of the Riemannian symmetric space $G/K$ is made.
But in our case:
\bd
\i[(a)] we are not on the Riemannian symmetric space $G/K$, but on
some other homogeneous space $G/\tilde{K}$, $\tilde{K}\subset K$, $\tilde{K}
\neq K$;
\i[(b)] as the reader can check, it is possible to project on the
symmetric space  $G/K_{g_{0}}$, but this problem is not invariant by the action of
$G$ on $G/K_{g_{0}}$. In fact, we get a {\it singular Riemannian} problem over
$G/K_{g_{0}}$ (as we have shown in the paper \cite{BCG}).
\ed  
\erem
With this methods,
for the real-resonant problem we get in
a more natural  way the result of \cite{BCG}, and we give
much simpler expressions for geodesics and optimal controls. 
In the general-complex problem, besides finding explicitly
expression for optimal trajectories and  controls we prove  
 our main result:
\bt {\bf (Main Result)} For the three-level problem 
with 
complex controls,
optimality implies resonance. More precisely, controls
$\Omega_1(.),\Omega_2(.)$ are optimal if and only if they have the
following form:
\bqn 
\left\{\ba{l} 
\Omega_1(t)=\cos(t/\sqrt{3})e^{i[(E_2-E_1)t+\varphi_1]},\\
\Omega_2(t)=\sin(t/\sqrt{3})e^{i[(E_3-E_2)t+\varphi_2]}. \ea\right.
\eqnl{main-res}
where $\varphi_1,\varphi_2$ are two arbitrary phases. 
Here the final time $T$
is  fixed in such a way subriemannian geodesics are parameterized by
arclength, and it is given by $T=\frac{\sqrt{3}}{2}\pi$.
\llabel{t-main}
\et
This fact was  pointed out as an open question in \cite{BCG}.
In other words, {\it optimal trajectories for the real-resonant
problem are
optimal also for the general-complex one}.
\brem
The fact that optimality implies resonance was also proved in \cite{BCG}
for
the two-level case, reducing the problem to an isoperimetric problem. More
precisely the two-level problem, which is described by the Hamiltonian:
\bqn 
H=\left(\ba{cc} E_1&\Omega(t)\\ \Omega^{\ast}(t)&E_2 \ea\right).
\eqnn
can be  reduced to a 3 dimensional contact
subriemannian problem, with a
special feature: it has a symmetry, transverse to the distribution. It is
a standard fact that such a subriemannian problem is  an
isoperimetric problem (in the sense of the calculus of variations) on the
quotient by the symmetry. 
In this case, the
fact that optimality implies resonance
is nothing but the
(trivial) solution of the classical isoperimetric
problem (or Dido problem) on the Riemannian sphere. We refer to \cite{BCG}
for details.
\llabel{r-iso}
\erem
\brem
Besides physical applications, the three-level problem is particularly
interesting since it has a
nontrivial geometric structure, but it is completely computable.
On the other side  the
two-level problem, described in Remark \ref{r-iso}
(although it is a $\k\oplus\p$ problem, with a nice geometric 
description) has a quite trivial solution.

Systems with more than 3 levels appear to be more difficult to treat.
The drift can be eliminated in the same way, but one gets a problem in
which the control distribution is determined by a strict subspace of $\p$ 
only.
In this case the Jurdjevic's formula fails to be true, 
and the integrability of the PMP is an open problem. 
See
Remarks 
\ref{r-nlevel2}  and 
\ref{r-appendix} (the last in Appendix B).  
In this case the Hamiltonian reads:
\bqn 
H&=& \left( \begin{array}{ccccc} E_{1} & \Omega
_{1}(t) & 0 & \cdots & 0 \\ \Omega_{1}^{\ast}(t) & E_{2} & \Omega _{2}(t)
& \ddots & \vdots \\ 0 & \Omega _{2}^{\ast}(t) & \ddots & \ddots & 0 \\
\vdots & \ddots & \ddots & E_{n-1} & \Omega _{n-1}(t) \\ 0 & \cdots & 0 &
\Omega _{n-1}^{\ast}(t) & E_{n} \end{array}
\right).
\llabel{ham-gen-n} 
\eqn 
\erem
\brem
For questions related to invariance 
under time
reparameterization of our optimal trajectories, 
 we refer to \cite{BCG}.
\erem
\subsection{A crucial fact}\llabel{s-crucial}
At this level, we want to point out an obvious, but important symmetry 
property. It is important for computations in the following, but also for 
practical applicability of the result.

The optimal control problem we consider (in $\C^3$ or in the complex 
sphere in $\C^3 $) is invariant under multiplication by a constant phase 
$e^{i\varphi}$. Hence, the problem of minimizing the energy from a fixed 
point in the ground eigenstate, to the third eigenstate, is the same as 
minimizing the energy to move from the full ground eigenstate, to the same 
target: actually, if $(\psi(t),\Omega(t))$ is an optimal pair 
trajectory-control, from 
$\psi_0$  in the ground eigenstate to the third one, 
$(e^{i\varphi} \psi(t),\Omega(t))$ is an optimal trajectory from 
$e^{i\varphi}\psi_0$. Moreover, as a consequence, whatever the initial 
condition $\psi_0$ in the ground eigenspace, the optimal control is 
\underline{the same}. This is very important for applicability. Although 
it 
is a time dependent problem, these considerations hold also for the 
real-resonant problem.  
\subsection{The problem downstairs and upstairs}
The problem of inducing a transition from the first to the third 
eigenstate, can be formulated, as
usual,
at the level of the wave function $\psi(t)$ but also at  the level of the
time evolution operator (the resolvent), 
 denoted here by $g(t)$: 
\bqn\psi(t)=g(t)\psi(0),~~~g(t)\in U(3),~~~g(0)=id.
\eqnl{temp-evol}
In the following we will call the optimal control problem for $\psi(t)$
and for $g(t)$ respectively {\it the "problem downstairs"} and {\it the 
"problem upstairs".}

The state-vector $\psi(t)$, solution of the time-dependent Schr\"odinger
equation $i\dot\psi=H\psi$, where $H$ is given by formula \r{ham-gen-3}
can be expanded in the canonical basis
of $\C^3$, formed by elements
$\varphi_1=(1,0,0),~\varphi_2=(0,1,0),~\varphi_3=(0,0,1)$:
$\psi(t)=c_1(t)\varphi_1+c_2(t)\varphi_2+c_3(t)\varphi_3$, with
$|c_1(t)|^2+|c_2(t)|^2+|c_3(t)|^2=1$. For $t<0$ and $t>T$,
$|c_i(t)|^2$ is the probability of measuring energy $E_i$.
Notice that, since $\Omega_j(t)=0$, for all $t<0$, $t>T$, $j=1,2$, we 
have:  
$$\frac d{dt}|c_i(t)|^2=0~~~\mbox{ for } 
t<0\mbox{ and }t>T.
$$   
At the level of the wave function, we formulate  the problem in 
following way. Assuming
$|c_1(t)|^2=1$ for $t<0$, we want to determine suitable control
functions $\Omega_j(.),$ $j=1,2$, such that $|c_3(t)|^2=1$ for time
$t>T$, requiring that they minimize
the cost (\ref{fluence-om}). Thus we have a control problem on the real
sphere $S^5\subset\C^3$ with initial point belonging to the circle
${\cal S}_\C^d$ defined by 
$|c_1|^2=1$ and target ${\cal T}_\C^d$ defined by $|c_3|^2=1$.
Equivalently, as we said in \ref{s-crucial}, the initial point $\psi(0)$ 
can be considered as free in ${\cal S}_\C^d$.
In the following the labels $(^d)$ and $(^u)$, indicate respectively 
downstairs and upstairs. 
Sources and targets upstairs (that will be called ${\cal S}_\C^u$ and ${\cal 
T}_\C^u$) will be computed in Section 
\ref{s-sources} after elimination  of the  drift. Why we use the subscript 
$(_\C)$, will be made clear in Remark \ref{r-notation}.

\subsection{Contents of the Paper}
The paper is organized as follows. 
In Section \ref{s-k+p}, we define precisely the $\k\oplus\p$ problem and  
in Section \ref{s-drift}, we discuss the elimination of the drift term. 
In Section \ref{s-formal} we formulate our problems in the $\k\oplus\p$
form
and we define  sources and the targets.
In Section \ref{s-geo} we compute optimal trajectories reaching the final
target and satisfying transversality conditions. The proof of Theorem 
\ref{t-main} 
follows. In Appendix A we point out an
interesting consequence of the Cartan decomposition, while in Appendix B
and C respectively, we explain the Jurdjevic's formalism and why abnormal
extremals are not optimal in a semi-simple $\k\oplus\p$ problem. 
\section{The $\k\oplus\p$ Problem\llabel{s-k+p}}
For sake of simplicity in the exposition, all over the paper, when we deal
about Lie groups and Lie algebras, we always consider that they are groups
and algebras of matrices.

Let $\l$ be a semi-simple Lie algebra and let us denote the Killing form 
by 
$Kil(.,.)$, $Kil(X,Y)=Tr(ad_X\circ ad_Y)$. In the following we recall what
we mean by a Cartan decomposition of $\l$.      

\bdeff 
A Cartan decomposition of a semi-simple Lie algebra
$\l$ is any
decomposition of the form:
\bqn
\l=\k\oplus\p, 
 \mbox{ where } [\k,\k]\subseteq\k,~~ [\p,\p]\subseteq\k,~~
[\k,\p]\subseteq\p.  
\eqnl{cartan-dec}
\llabel{d-k+p}
\edeff

\vspace{-1cm}
\brem
Since $\l$ is semi-simple then relations \r{cartan-dec} implies 
$[\p,\p]=\k,~~
[\k,\p]=\p$ (see Appendix A for the proof). This fact will be  crucial for
the elimination of abnormal extremals. 
\erem
\bdeff 
The right-invariant $\k\oplus\p$ control problem on a \underline{compact}
semi-simple Lie group
$G$ is the
subriemannian problem with right-invariant distribution induced by $\p$ 
and cost:
$$
\int_{0}^{T}<\dot g~g^{-1},\dot g ~g^{-1}> dt,
$$ 
where $<~,~>:=-\al~Kil\big|_p(~,~)$, $\al>0$. 
\llabel{d-k+p-problem}
\edeff 
The constant $\al$ is clearly not relevant. It will be used just to obtain
good normalizations.

In the following we will be interested in $\k\oplus\p$ problems on $so(3)$
and
$su(3)$, for which we have respectively $Tr(ad_X\circ ad_Y)=Tr(X~Y)$ and
$Tr(ad_X\circ ad_Y)=5~Tr(X~Y)$.  Then, in order to
get for both Lie 
algebras the useful relation:
\bqn
<X,Y>=-\frac12 Tr(X~Y),~~X,Y\in\p,
\eqnl{metric}
we must set $\al=1/2$ for $so(3)$ and $\al=1/10$
for $su(3)$.

Let $\{X_j\}$ be an orthonormal (right-invariant) frame for the
subspace $\p$ 
of $\l$, with
respect to the
metric defined in Definition \ref{d-k+p-problem}. 
Then the $\k\oplus\p$ problem
reads:
\bqn
\dot g=\left(\sum_j u_jX_j\right)g,~~~~
min \int_{0}^{T}\sum_ju_j^2~dt,~~~~u_j\in\R.
\eqnl{control-p}
From relations \r{cartan-dec}, one gets that the so called {\it "Goh 
condition"} is
never
satisfied (see Appendix C). As a consequence:
\bp
In the problem defined in Definition \ref{d-k+p-problem}, every 
strict abnormal
extremal (if any) is not optimal. 
\ep
\brem
The $\k\oplus\p$ problem can be  also stated  in the case where $\l$ is of
noncompact type. In this case, to have a positive definite 
metric, we
have to require  
that $\k$ is a maximal compact subalgebra of
$\l$,  and we must define 
 $<~,~>:=+\al~Kil\big|_p(~,~)$, $\al>0$.  
\llabel{r-maximalcompact}
\erem
\section{Elimination of the Drift Term \llabel{s-drift}}
In this section, we show how to eliminate the drift term from a $n$-level 
system of the form \r{ham-gen-n} with general-complex and real-resonant 
controls. 
For $\Omega\in\C$, let us denote by $M_j(\Omega)$ and $N_j(\Omega)$ the 
$n\times n$ matrices:
\bqn
\left\{\ba{l}
M_j(\Omega)_{k,l}=\delta_{j,k}\delta_{j+1,l}\Omega+\delta_{j+1,k} 
\delta_{j,l}\Omega^\ast\\
N_j(\Omega)_{k,l}=\delta_{j,k}\delta_{j+1,l}\Omega-\delta_{j+1,k}, 
\delta_{j,l}\Omega^\ast
\ea\right.,~~j=1,...,n-1,
\eqnl{belin}
where $\delta$ is the Kronecker symbol: $\delta_{i,j}=1$ if $i=j$,  
$\delta_{i,j}=0$ if $i\neq j$.
Let $\Delta=diag(E_1,...,E_n)$, $\omega_j=E_{j+1}-E_j$, $j=1,...,n-1$. We 
will consider 
successively the \underline{general complex} problem (in dimension $n$):
\bqn
i\dot\psi=H\psi,~~~H=\Delta+\sum_{j=1}^{n-1}M_j(\Omega_j),
\mbox{ where $\Omega_j\in\C$,}
\eqnn
(this is nothing but another notation for the matrix \r{ham-gen-n}),
and the \underline{real-resonant} problem:
\bqn
H=\Delta+\sum_{j=1}^{n-1}M_j(e^{i(\omega_j 
t+\alpha_j)}u_j),~~~u_j,\al_j\in\R.
\eqnn
In both problems, $\psi$ lies in the complex sphere in $\C^n$, and we 
want to connect the source ${\cal S}_\C^d=\{(e^{i\varphi},0,...,0)\}$ to the 
target ${\cal T}_\C^d=\{(0,...,0,e^{i\varphi})\}$, by minimizing:
\bqn
J=\int_0^T\sum_{j=1}^{n-1}|\Omega_j|^2dt,
\mbox{ which in the real-resonant case is  }
J=\int_0^T\sum_{j=1}^{n-1}u_j^2dt.
\eqnn
In both cases, we first make the change of variable $\psi=e^{-i\Delta 
t}\Lambda$ (interaction representation), to get (here  $Ad_\phi B:=\phi 
B\phi^{-1}$):
\bqn
i\dot\Lambda=\sum_{j=1}^{n-1}\big(Ad_{e^{i\Delta 
t}}M_j(\Omega_j)\big)\Lambda=\sum_{j=1}^{n-1}M_j\big(e^{-i\omega_j t 
}\Omega_j)\Lambda.
\eqnn
Let us stress that the source $\Sdc$ and the target $\Tdc$ are preserved by this first 
change of 
coordinates.
\subsection{The general-complex case}
In that case, we make  the time-dependent {\it gauge transformation} 
(i.e. cost preserving change of controls):
\bqn
e^{-i\omega_j t}\Omega_j=i\tilde \Omega_j.
\eqnn
Hence our problem become (after the change of notation 
$\Lambda\to\psi$, $\tilde \Omega_j \to u_j$):
\bqn
\left\{\ba{ll}
 \mbox{\textbf{a.}} \  & \min\int_0^T\sum_{j=1}^{n-1}|u_j|^2dt,~~x(0)\in
  {\cal S}_\C^d,~~x(T)\in 
{\cal T}_\C^d, \\
 \mbox{\textbf{b.}} \  & \dot \psi=\sum_{j=1}^{n-1}N_j(u_j)\psi,~~ 
u_j(t)\in\C. 
\ea
\right.
\eqnl{ham-gen-u}
Notice that the matrices $N_j(1), N_j(i)$ generate $su(3)$ as a Lie algebra.
For $n=3$, the \sceq (\ref{ham-gen-u}\textbf{b}) writes, in matrix form:
\bqn
\dot\psi=\tilde H_\C\psi\mbox{ ~~~where~~~ }\tilde H_\C:=
\left(\ba{ccc} 0& u_1(t) &0\\ - u_1^\ast(t) &0 &u_2(t) \\
0&- u_2^\ast(t) &
0
\ea\right).
\eqnl{ham-c-i}   
The cost and the relation between  controls before and after elimination 
of the drift are:
\bqn
&&J=
\int_{0}^{T}\left(|u_1(t)|^2+|u_2(t)|^2\right)~dt,\llabel{fluence-u}\\
&&\left\{\ba{l}
\Omega_1(t)=u_1(t)e^{i[(E_2-E_1)t+\pi/2]},\\
\Omega_2(t)=u_2(t)e^{i[(E_3-E_2)t+\pi/2]},
\ea\right.
\llabel{ND2} 
\eqn

\subsection{The real-resonant case}
In this case, since $\Omega_j=u_j e^{i(\omega_j t+\al_j)}$, we have:
\bqn
i\dot\Lambda=\sum_{j=1}^{n-1} 
M_j\big(e^{i\al_j}u_j\big)\Lambda,~~~u_j\in\R.
\eqn    
We make another diagonal, linear change of coordinates:
\bqn
\Lambda=e^{i L}\phi,~~L=diag(\lam_1,...\lam_n),\mbox{ for 
}\lam_1,...,\lam_n\in\R.
\eqnn
This gives:
\bqn
i\dot\phi=\sum_{j=1}^{n-1}M_j\big(e^{i(\al_j+\lam_{j+1}-\lam_j)}u_j\big)\phi.
\eqnn
Choosing the $\lam_j$'s for $e^{i(\al_j+\lam_{j+1}-\lam_j)}=i$, we get:
\bqn
\dot\phi=\sum_{j=1}^{n-1}N_j(u_j)\phi,~~u_j(t)\in\R.
\eqnl{6jp}
The source and the target are also preserved by this change of 
coordinates. 
Notice that the matrices $N_j(1), j=1 \ldots n-1$ in \r{6jp},
generate  $so(n)$ as a Lie algebra in its presentation by skew-symmetric
 matrices. This means that the orbit of the 
system \r{6jp} through the points $(\pm1,0,...,0)$ is the real sphere 
$S^{n-1}$. 
In other words, the action  
on $\C^n$ of the subgroup $SO(n)\subset SU(n)$, restricts to the 
reals.
Hence (by multiplication on the right by 
$e^{i\varphi}$), the orbit through the points  $(\pm 
e^{i\varphi},0,...,0)$ 
is 
the set $S^{n-1}e^{i\varphi}$. Therefore, (after the change of notation
$\phi\to\psi$) the 
real-resonant problem reduces to 
the problem over $S^{n-1}$ (see sections \ref{s-controllability} 
for more details):
\bqn
\left\{\ba{l}
\min\int_0^T\sum_{j=1}^{n-1}u_j^2dt,~~x(0)\in \{(\pm1,0,...,0)\},~~x(T)\in 
\{(0,...,0,\pm1)\},\\
\dot \psi=\sum_{j=1}^{n-1}N_j(u_j)\psi,~~ u_j(t)\in\R.
\ea
\right.
\eqnn
For $n=3$ we get for the \sceq in matrix form: 
\bqn
\dot\psi=\tilde H_\R\psi\mbox{ ~~~where~~~ }\tilde H_\R:
=\left(\ba{ccc} 0& u_1(t) &0\\ -  u_1(t) &0 & u_2(t) \\
0&-u_2(t) &
0
\ea\right),
\eqnl{ham-r-i}
The cost is given again by formula \r{fluence-u} and 
the relation between  controls before and after elimination 
of the drift is:
\bqn 
&&\Omega_j(t)=u_j(t)~e^{i (\omega_j t+\al_j)},~~\omega_j=E_{j+1}-E_j,\nn\\
&&u_j(.):\R\to\R,~~\al_j\in[-\pi,\pi],~~j=1,2.
\eqnl{ND1}

\brem
In the following, besides to the labels 
$(^d)$ and $(^u)$ that indicate respectively 
downstairs and upstairs, we will use the labels $(_\C)$ and $(_\R)$  to 
indicate respectively the general-complex problem and the real-resonant 
one. When these labels are dropped in a formula, we mean that it is valid 
for both the real-resonant and the general-complex problem. With this 
notation:
\bqn
\Sdc=\{(e^{i\varphi},0,0)\},~~~\Tdc=\{(0,0,e^{i\varphi})\},\nn\\
\Sdr=\{(\pm1,0,0)\},~~~\Tdr=\{(0,0,\pm1)\}.\nn
\eqn
\llabel{r-notation}
\erem
\section{The $\k\oplus\p$ Form of the Problem, 
Controllability, Sources and Targets\llabel{s-formal}}
\subsection{Controllability}\llabel{s-controllability}
First notice that, after elimination of the drift, 
the problem upstairs is in $SU(3)$ and not in 
$U(3)$ (the center of $U(3)$ is eliminated). Moreover, as we already 
observed, in the real resonant case, 
starting from $(1,0,0)$,
we have the standard action of $SO(3)$ on $\R^3$. From the fact that the
matrices $\tilde H_\R$
generate 
$so(3)$ as a Lie algebra, it follows easily that the orbit through 
$(1,0,0)$ is the real sphere
$S^2$ 
of equation: 
\bqn
c_1^2+c_2^2+c_3^2=1,~~~c_j\in\R,~~~j=1,2,3.
\eqnl{orbitasferica}

It is equivalent to say that  the real-resonant problem is not 
controllable on $S^5$. 
This fact is explained in details in \cite{BCG}. Moreover in  \cite{BCG}
it is proved that 
on each of these spheres $S^2$ the control problem reduces to a
singular-Riemannian problem. The ``relevant locus'',
which is the union of all the orbits (translations of $S^2$) passing through the
eigenstate number 1, has an interesting non trivial geometric description. It is the
only non orientable sphere-bundle over $S^1$. 

After elimination of the drift, the real-resonant problem upstairs is on 
$SO(3)$, and the general-complex problem upstairs is in $SU(3)$. We have:
\bp
After elimination of the drift, the real resonant problem upstairs is 
completely controllable on $G_\R:=SO(3)\subset U(3)$, and the 
general-complex 
problem is completely controllable on $G_\C:=SU(3)\subset U(3)$. As a 
consequence, by 
transitivity of the action on the real and complex spheres, the 
corresponding problems downstairs are controllable on $S^2$ and $S^5$.
\llabel{p-controll}
\ep
\proof
This is a consequence of the general theorems of controllability of 
left-invariant control systems on compact groups: the {\it Lie-rank 
condition} is necessary and sufficient for controllability (see 
\cite{jurd-suss}). \quadp\\
\brem In fact, before the elimination of the drift, the general 
complex 
problem upstairs is in $U(3)$, due to the nonzero trace of $\Delta$. For 
the same reason as in Proposition \ref{p-controll}, it is completely 
controllable.
\erem 

\subsection{The problems in the $\k\oplus\p$ form}
For the lifted problem,
from $\dot\psi=\tilde H\psi$, using \r{temp-evol} one gets:
\bqn
\dot g =\tilde H g.
\eqnl{gHg}
where $\tilde H_\R$ generates $so(3)$ and $\tilde H_\C$ 
generates $su(3)$ (as Lie algebras, for distinct values of the controls). 
For both Hamiltonians, the problem of minimizing the cost \r{fluence-u} is
a $\k\oplus\p$ problem.\\ In fact in the real-resonant case equation \r{gHg}
can
be
written as:
\bqn
\dot g =(u_1 X_1+u_2 X_2) g,
\eqnn
where:
\bqn
X_1=\left(\ba{ccc} 0&1 &0\\  -1 &0 &0 \\ 0&0 &0\ea\right),~~
X_2=\left(\ba{ccc} 0&0&0\\ 0 &0 &1 \\ 0&-1 &0\ea\right).
\eqnl{x1x2}
Setting: 
\bqn
X_3=\left(\ba{ccc} 0& 0&1\\  0 &0 &0 \\ -1&0
&0\ea\right),~~~~\p:\mbox{=span(}\{X_1,X_2\}\mbox{)},~~\mbox{ and 
}~~\k=\mbox{span}(\{X_3\}),
\eqnn 
one gets relations \r{cartan-dec} (with $\l=so(3)$). 
Moreover the distribution is right-invariant and the frame \r{x1x2}  
is orthonormal
for the metric \r{metric}.

In the case of $su(3)$ we have: 
\bqn
\dot g =(u_1 X_1+u_2 X_2+u_3 Y_1+u_4 Y_2) g,~~~u_j\in\R,
\eqnn
where $X_1$ and $X_2$ are given by formula \r{x1x2} and $Y_1$, $Y_2$
by:
\bqn
Y_1=\left(\ba{ccc} 0&i &0\\  i &0 &0 \\ 0&0 &0\ea\right),~~
Y_2=\left(\ba{ccc} 0&0&0\\ 0 &0 &i \\ 0&i &0\ea\right).
\eqnn
One can easily check that relations \r{cartan-dec} are
verified with $\p:=$span$(\{X_1,X_2,Y_1,Y_2\})$ and 
$\k:=$span$(\{Z_1,Z_2,Z_3,Z_4\})$, where:   
\bqn
Z_1=\left(\ba{ccc} 0&0 &1\\  0 &0 &0 \\ -1&0 &0\ea\right),~~
Z_2=\left(\ba{ccc} 0&0 &i\\  0 &0 &0 \\ i&0 &0\ea\right),~~
Z_3=\left(\ba{ccc} i&0 &0\\  0 &-i &0 \\ 0&0 &0\ea\right),~~
Z_4=\left(\ba{ccc} 0&0 &0\\  0 &i &0 \\ 0&0 &-i\ea\right),
\eqnn
are the  remaining 4 generators of $su(3)$.
Again the distribution is right-invariant and the frame 
$\{X_1,X_2,Y_1,Y_2\}$ orthonormal
for the metric \r{metric}.
\brem $\bullet$
 (general complex case) $\k$ is a subalgebra of $su(3)$ containing a 
maximal Abelian 
subalgebra
of  $su(3)$, which is generated by  $Z_3$ and $Z_4$.  From \cite{borel}
(see also \cite{lie-for-p})  it
follows that $\k$ must be a Borel subalgebra.  In this case, it is
$u(1)\times su(2)=u(2)$. 
\llabel{r-B}
\erem
\noindent
Let us denote by
$S(U(1)\times U(2))$ (resp. $S(Z_2 \times O(2))$ , the groups 
of matrices of the form: 
$$
B=\left(\mbox{  
\begin{tabular}{c|c} 
$\epsilon$&0 \\\hline
~&~\\
0&$~U_{\epsilon}~$ \\
~&~\\
\end{tabular} 
}
\right),\mbox{ where }\det(B)=1,
$$
with $\epsilon \in U(1)$, (resp. $\epsilon \in Z_2=\{-1,1\}$), and $U_{\epsilon} \in U(2)$ (resp. $U_{\epsilon} \in O(2)$), 
$S(U(1)\times U(2)) \approx U(2)$, $S(Z_2\times O(2)) \approx O(2)$.\\
\\\\
Set: 
\bqn
g_0:=\left(\ba{ccc}0&1&0\\0&0&1\\1&0&0\ea\right),
\llabel{g_0-15avril}
\eqn  
$g_0 \in SO(3) \subset SU(3)$,
 ${g_0}^{2}={g_0}^{-1}$. 
In the real-resonant case, we will denote by $K^{\R}$ the subgroup of $G^{\R}=SO(3)$, 
conjugate by ${g_0}^{-1}$ to $S(Z_2 \times O(2))\approx O(2)$, with Lie algebra
${\k}^{\R}={\{X_3\}}_{LA}$.\\
In the general-complex case, we denote by $K^{\C}$ the subgroup of 
$G^{\C}=SU(3)$ conjugate by
${g_0}^{-1}$ to $S(U(1) \times U(2)) \approx U(2)$ with Lie algebra
${\k}^{\C}={\{Z_1,Z_2,Z_3,Z_4\}}_{LA}$.
\bi
\i In the general-complex case, the corresponding Riemannian symmetric space $G/K$ is
$SU(3)/ S(U(2)\times U(1))\sim\C P^2$ that has rank $1$ (see
Helgason \cite{helgason}, pp.
518), and the Cartan subalgebra of $\p$ has dimension $1$ only, as one can check
easily. 
\i In the real-resonant case, it is $\R P^2$, the 2-dimensional projective 
space.
\ei
\brem
Notice that in the real-resonant problem  the distribution is a  
contact distribution. In this case 
it is a standard fact (see for instance
\cite{dido,andrei1,bellaiche,chakir1,gromov}) that there are no
abnormal extremals.  On the other hand, the distribution of the general
complex problem is not a contact distribution and abnormal extremals do
exist. Anyway they are not optimal, due to the fact that the  {\it "Goh
condition"} (see Appendix C) fails to hold. 
\erem
\brem {\bf (n-level case)} 
For $n\geq4$ the Hamiltonian given in formula \r{ham-gen-u} generates 
$so(n)$ or  
$su(n)$ respectively 
with real or complex controls .
Anyway it 
never gives rise to a 
$\k\oplus\p$ problem since the distribution is only a strict subspace of $\p$.
As explained in Remark \ref{r-appendix} this fact causes the failure of  
the proof of
the integrability of the PMP.  \underline{No results in this paper 
about optimal trajectories can be easily generalized to 4 or more levels}.
For 4 or more levels
the integrability of the PMP, and the fact that  optimality implies
resonance, are open questions.
\llabel{r-nlevel2}
\erem
 \subsection{Sources and Targets}
\llabel{s-sources}
In this section, we describe  sources  and  targets
for the real-resonant and the general-complex problems upstairs.

As we already said, from section \ref{s-crucial}, sources and targets 
downstairs, for the general-complex problem and the real-resonant one,  
are respectively the circles
${\cal S}_\C^d$, 
${\cal T}_\C^d$, and the sets
$\Sdr=\{(\pm1,0,0)\}$, 
$\Tdr=\{(0,0,\pm1)\}$.

In both cases, we decide that the canonical 
projection $\pi:G\to G/\tilde{K}$ 
maps
$\tilde{K}$ to the point $(1,0,0)$ in the complex and real sphere in $\C^3$ (resp. 
$\R^3$). Here $\tilde{K}_\C=SU(2)$ (general complex case) and   $\tilde{K}_\R=SO(2)$ 
(real resonant case). Let ${\cal S}^u=\pi^{-1}({\cal S}^d)$. Then 
$\Suc=S(U(1)\times 
U(2)),$ $\Sur=S(Z_2\times O(2))$. \\
The element $g_0 \in SO(3)$ defined in the previous section maps ${\cal S}_\C^d$ to
 ${\cal T}_\C^d$. Then 
    ${\cal T}^u:=\pi^{-1}({\cal T}^d)=g_0{\cal S}^u$ in both cases as it 
is easy to check. 


 It is clear that a trajectory $g(t)$ of the system upstairs:
 $$\begin{array}{lr}
   \dot g = \tilde{H}g, & t \in [0, T]
   \end{array} $$
 such that $g(0) \in {\cal S}^{u}$ , $g(T) \in {\cal T}^{u}$, maps into a 
trajectory 
 $g(t)\tilde{K}$ of the system on the sphere,$$
 \begin{array}{lcr}
 	\dot{g}\tilde{K}=\tilde{H}g\tilde{K} & g(0)\tilde{K}\in {\cal 
S}^d,&
	g(T)\tilde{K} \in {\cal T}^d 
 \end{array}	
	$$
with the same control, hence the same cost.	
   Conversely, if $x(t)$, $t \in [0, T]$,
 is a trajectory on the sphere, of the system 
   $\dot x = \tilde{H}x$, that maps the point $(\varepsilon, 0, 0), 
  |\varepsilon|=1$
   to the point $(0,0,\varepsilon'), |\varepsilon'|=1$.
   Then, the corresponding fundamental matrix solution $g(t)$ satisfies, 
   for all $s \in {\cal S}^u$, $g(0) s =s \in {\cal S}^u$ :  
   $$g(T).s.(\varepsilon, 0,0) =(0, 0,\varepsilon')=g_{0}.(\varepsilon',0,0)$$
   Therefore, $g_{0}^{-1}{g(T)}.s.(\varepsilon,0,0)=(\varepsilon',0,0)$ and
   $g_{0}^{-1}.{g(T)}.s \in {\cal S}^u$, $g(T) \in g_{0}{\cal S}^u={\cal 
T}^u$.
 This shows that, for all $s \in {\cal S}^u$, the lifted solution upstairs
 satisfies $g(0).s \in {\cal S}^u$, $g(T).s \in {\cal T}^u$.


 Hence we get 
the following table:
\footnotesize
\begin{center}
\begin{tabular}{|c||l|l|} \hline
PROBLEM&SOURCE&TARGET \\\hline
&&\\
real-resonant   &
$\Sur:=\left\{\left(\mbox{  
\begin{tabular}{c|c} 
$Z_2$&0 \\\hline
~&~\\
0&$~O(2)~$ \\
~&~\\
\end{tabular} 
}
\right)
\in SO(3)\right\}=S(Z_2\times O(2))$&
$\Tur:=g_0~\Sur=g_0 S(Z_2\times O(2))$
\\
&&\\
\hline
&&\\
general-complex & 
$\Suc:=\left\{\left(\mbox{  
\begin{tabular}{c|c} 
$U(1)$&0 \\\hline
~&~\\
0&$~U(2)~$ \\
~&~\\
\end{tabular} 
}
\right)
\in SU(3)\right\}=S(U(1)\times U(2))$&
$\Tuc:=g_0\Suc=g_0S(U(1)\times U(2))$
\\
&&\\
\hline
\end{tabular} 
\end{center}
\normalsize

\section{Expression for Geodesics, Controls and Transversality
Conditions\llabel{s-geo}}
\subsection{Preliminaries}
As explained in Appendices B and C, candidates optimal trajectories 
 for the $\k\oplus\p$ problem are only ``normal geodesics'' and are given
by
the following formula: 
\bqn 
g(t)=e^{-A_k t}e^{(A_k+A_p)t}g(0),
\eqnl{jurdjevic} 
where $g(0)$ is the starting point belonging to the
source. Set $M(0)=M_p(0)+M_k(0)=A=A_k+A_p\in\l$ for the initial value of the
covector $P(0)=d_{R_{g}}^{\ast}p_{g(0)}=\langle M(0),\cdot \rangle$, where
$\langle \cdot,\cdot \rangle$ is given in (\ref{metric}), and with the
notations of Appendix B. Here, $A_p=M_p(0)$, $A_k=M_k(0)$ are the orthogonal
projections of $M(0)=A$ on $\p$ and $\k$ respectively.
Geodesics are parametrized by arclength iff: 
\bqn
\langle A_p,A_p \rangle =1 
\eqnl{apap}
\subsection{Transversality conditions} 
  In both cases, our source is
  $K_{g_{0}} = g_{0} K {g_{0}}^{-1}$ where $g_0$ has been defined in 
formula
(\ref{g_0-15avril})
  and $K_{g_{0}}^{\C}=S(U(1)\times U(2))$,
   $K_{g_{0}} ^{\R}=S(Z(2)\times O(2))$. Let us notice the two following
   facts:\\\\
   \underline{\textbf{Facts:}}
   \begin{enumerate}
   \item Transversality conditions at the source may be required at the 
identity only.
   \item Transversality conditions at the source imply transversality
   conditions at the target.
   \end{enumerate}
\textbf {Proof} : Point 1 comes from right-invariance and the fact that 
the source is a subgroup. Point 2 comes from the
following lemma \ref{lemmeAd}. \quadp\\\\
Let $(g(t),\Omega(t))$, $t\in [0, T]$ be a normal extremal
	corresponding to the covector $p_{g(t)}$, such that $g(0)=Id$.
\bl $p_{Id}(Ad{g_0}\k)=0$ implies $p_{g(t)}(dL_{g(t)}Ad_{g_0}\k)=0$, $\forall
t \in [0, T]$.
\llabel{lemmeAd}
\el
\textbf {Proof} : (with the notations of Appendix B). Set
$
I=p_{g(t)}(dL_{g(t)}Ad_{g_0}\k)=dR_{{g(t)}^{-1}}^{\ast}P(t)
(dL_{g(t)}Ad_{g_0}\k).
$
Then $I=P(t)(Ad_{g(t)g_{0}}\k)=\langle M(t), Ad_{g(t)g_0}\k \rangle.$
But: 
\bqn
M(t)&=&M_p(t)+M_k(t)=M_p(t)+M_k(0),\nn\\
M(t)&=&e^{-M_k(0)t}M_p(0)e^{M_k(0)t}+ M_k(0)=
e^{-M_k(0)t}(M_p(0)+M_k(0))e^{M_k(0)t}=
Ad_{e^{-M_k(0)t}}M(0),\nn
\eqn
and $g(t)=e^{-M_k(0)t}e^{M(0)t}$. Then: 
$$I=\langle Ad_{e^{-M_k(0)t}}M(0), 
Ad_{e^{-{M_k}(0)t}}Ad_{e^{M(0)t}}Ad_{g_0}\k\rangle,$$
and since the killing form is $Ad_G$ invariant, $I=\langle M(0),$ 
$Ad_{e^{M(0)t}}Ad_{g_0}\k\rangle$. Hence, for the same reason:
\bqn
I =\langle{Ad_{e^{-M(0)t}}}M(0),Ad_{g_0}\k\rangle= \langle 
M(0),Ad_{g_0}\k\rangle = P(0)(Ad_{g_0}\k)= 
p_{Id}(Ad_{g_0}\k)=0.\nn
\eqn
\quadp\\
Similarly to $N_1$, $N_2$, (cfr. formula \r{belin}) 
let us define $N_{1,3}$ 
by:
\bqn 
N_{1,3}(a_3 e^{i {\theta}_3})=\left(\ba{ccc} 0&0&a_3e^{i {\theta}_3}\\ 0&0&0\\
-a_3 e^{-i {\theta}_3}&0&0 \ea\right).
\eqn
Let us set, in the real-resonant case $A_p=a_1 N_1(1)+ a_2 N_2(1)$ and $A_k=a_3 N_{1,3}(1)$.
In the general complex case, set 
$A_p=N_1(a_1 e^{i {\theta}_1})+N_2(a_2 e^{i {\theta}_2})$,
 $A_k=a_4 Z_3+a_5 Z_4+N_{1,3}(a_3 e^{i {\theta}_3})$.
Here, $a_i \in \R$ and ${\theta}_i \in [-\pi , \pi]$.

\subsection{The real-resonant case (G=SO(3))} 
\bp For the real-resonant problem, the transversality conditions
$Kil(A,T_{id}\Sur)=0$ implies $a_2=0$. \llabel{p-tc1} \ep
\proof We have: 
\bqn
Ad_{g_0}{\k}^{\R}=T_{id}\Sur:=\left\{\left(\mbox{ \begin{tabular}{c|cc} $0$&$0$&$0$
\\\hline $0$&$0$&$-\beta$\\ $0$&$\beta$&$0$ \end{tabular} } \right),
~~~\beta\in\R\right\}. 
\eqnn 
Then equation $Kil(A,T_{id}\Sur)=0$ is satisfied for
every $\beta\in\R$ if and only if $a_2=0$. \quadp \\
From
Proposition \ref{p-tc1}, and condition \r{apap}, one gets the covectors to 
be used in formula
\r{jurdjevic}:  \bqn A^\pm= \left(\ba{ccc} 0 &\pm1&a_3\\ \mp1&0 &0\\
-a_3&0&0 \ea\right),~ \eqnl{A1}
\bp The
geodesics \r{jurdjevic}, for which $g(0)=Id$, with $A$ given by formula 
\r{A1}, reach the
target $\Tur$ for the smallest time (arclength) $|t|$, if and only if $a_3=\pm1/\sqrt{3}$.
Moreover, the 4
geodesics
(corresponding $A^\pm$ and to the signs $\pm$ in $a_3$) have the same
length and reach the target at time:  $$T=\frac{\sqrt{3}}{2}\pi.$$
\llabel{p-final1} \ep
\proof
Computing $g(t)=e^{-A_k
t}e^{(A_k+A_p)t}$, with $A$
given by formula \r{A1},
 and recalling  that: 
$$\psi(t)=g(t)\psi(0)=g(t)\left(\ba{c}1\\0\\0\ea\right),$$
one gets for the square of the third component of the wave function:
\bqn
\big(c_3(t)\big)^2=
      {\frac{{{\left( \cos (t\,{a_3})\,
          \sin (t\,
            {\sqrt{1 + {{{a_3}}^2}}}
            )\,{{{a}}_3}\,
          {\sqrt{1 + {{{a_3}}^2}}} -
          \cos (t\,
            {\sqrt{1 + {{{a_3}}^2}}}
            )\,\sin (t\,{a_3})\,
          \left( 1 +
           {{{a_3}}^2} \right)
          \right) }^2}}{{{\left(
         1 + {{{a_3}}^2} \right) }^2
      }}}
\eqnl{c3el2}
 By lemma (2) in Appendix D, we get the result. \quadp \\\\
{\bf Explicit expressions for the wave function and for optimal
controls}\\
Let us fix for instance the sign $-$ in \r{A1} and $a_3=+1/\sqrt{3}$.
The expressions of the three components of the wave function are:
\bqn
\left\{\ba{l}
c_1(t)={{\cos ({\frac{t}{{\sqrt{3}}}})}^3}\\
c_2(t)=\frac{\sqrt{3}}{2}\,
       \sin ( \frac{2\,t}
          {\sqrt{3}} )   \\
c_3(t)=-{{\sin ({\frac{t}{{\sqrt{3}}}})}^3}
\ea\right.
\eqnl{RC1} 
 Let us stress that this curve is not a circle on $S^2$.\\
 \\
Controls can be obtained with the following expressions:
\bqn
 u_1=(\dot g g^{-1})_{1,2},~~ u_2=(\dot g g^{-1})_{2,3}
 \eqnl{labelX}\\
 We get:
\bqn
\left\{\ba{l}
u_1(t)=
-\cos ({\frac{t}{{\sqrt{3}}}})\\
u_2(t)=
\sin ({\frac{t}{{\sqrt{3}}}})
\ea\right.
\eqnl{RU1} 
The situation is depicted on Figure 1.
\ppot{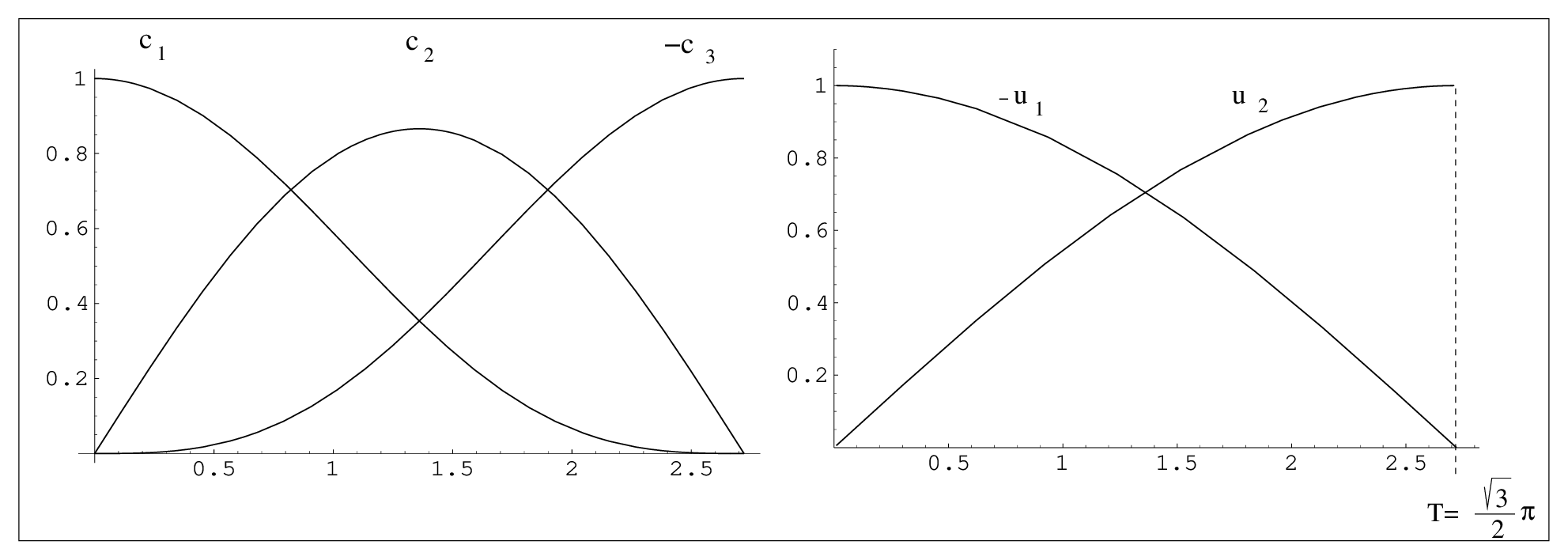}{Figure 1}{6}{15}
The result given by formulas \r{RC1} and \r{RU1} is the same result as the one of
\cite{BCG}, but it has a simpler form. Using expression \r{ND1} (resonance
hypothesis), we get for the external fields:
\bqn
\left\{\ba{l}
\Omega_1(t)=-\cos(t/\sqrt{3})e^{i(\omega_1t+\alpha_1)},\\
\Omega_2(t)=\sin(t/\sqrt{3})e^{i(\omega_2t+\alpha_2)}.
\ea\right.
\eqnl{OM1}
Notice that the phases $\alpha_1,~\alpha_2$ are arbitrary.
\subsection{The general-complex case: G=SU(3)}

\bp
For the general-complex problem, the transversality conditions 
$Kil(A,T_{id}\Suc)=0$ 
implies 
$a_2=a_4=a_5=0$.
\llabel{p-tc2}
\ep
\proof
We have:
\bqn
Ad_{g_0}\k^{\C}=T_{id}\Suc:=\left\{\left(\mbox{  
\begin{tabular}{c|cc} 
$i\al_1$&$0$&$0$ \\\hline
$0$&$i(\al_2-\al_1)$&$ \beta_1+i\beta_2$\\
$0$&$-\beta_1+i\beta_2$&$-i\al_2$
\end{tabular} 
}
\right), ~~~\al_1,\al_2,\beta_1,\beta_2\in\R\right\}
\eqnn	
Then equation $Kil(A,T_{id}\Suc)=0$, is satisfied for every
$\al_1,\al_2,\beta_1,\beta_2\in\R$ if and only
if 
$a_2=a_4=a_5=0$.
\quadp \\ \\
The covector to be used in formula
\r{jurdjevic} is then:
\bqn
A^{(\theta_1,\theta_3)}=
\left(\ba{ccc} 
0   &e^{i\theta_1}&a_3e^{i\theta_3}\\ 
-e^{-i\theta_1}&0   &0\\ 
-a_3e^{-i\theta_3}&0&0
\ea\right).
\eqnl{A2}
\bp 
The
geodesics \r{jurdjevic}, with $A$ given by formula \r{A2} (for which $g(0)=Id$), reach the
target $\Tuc$ for the smallest time (arclength) $|t|$, if and only of $a_3=\pm  1/\sqrt{3}$. Moreover 
all the geodesics of the two parameter family corresponding
to
$\theta_1,\theta_3\in[-\pi,\pi]$, have the same length:  
$$T=\frac{\sqrt{3}}{2}\pi.$$\llabel{p-final2}
\ep 
\proof 
The explicit expression for $|c_3|^2$ is given by the right-hand side of 
formula \r{c3el2}. The conclusion follows as in  the proof of
Proposition \r{p-final1}. \\ 
\\
{\bf Explicit expressions for the wave function and for optimal
controls}\\
The expressions of the three components of the wave function and of 
optimal controls are: 
\bqn 
\left\{\ba{l}
c_1(t)={{\cos ({\frac{t}{{\sqrt{3}}}})}^3}\\
c_2(t)=-\frac{\sqrt{3}}{2}\,
       \sin ( \frac{2\,t}
          {\sqrt{3}} )e^{-i\theta_1}   \\
c_3(t)=-{{\sin ({\frac{t}{{\sqrt{3}}}})}^3} e^{-i\theta_3}, 
\ea\right.
\llabel{RC2}\\ 
\left\{\ba{l} 
u_1(t)=\cos(t/\sqrt{3})e^{i\theta_1}   \\
u_2(t)=-\sin(t/\sqrt{3})e^{i(\theta_3-\theta_1)},
\ea\right.
\eqnl{RU2} 
\brem 
Again, notice that none of these trajectories is a circle on the
corresponding (translation of) sphere $S^2$.
Notice that all the geodesics of the family described by
Proposition (\ref{p-final2}) have the same length as the 4 geodesics
described by Proposition
\ref{p-final1}. This proves that the use of the complex Hamiltonian
\r{ham-c-i} (instead of the real one \r{ham-r-i}) does not allow to
reduce
the cost \r{fluence-u}, and this proves Theorem \ref{t-main} in the 
introduction.
\erem
 Formulas \r{RC1}, \r{RU1} can be obtained from
formulas \r{RC2}, \r{RU2} setting $\theta_1=\pi,~~\theta_3=0$. \\

For the general-complex problem, using expressions \r{RU2} in \r{ND2} one
gets for the external fields:
\bqn
\left\{\ba{l}
\Omega_1(t)=\cos(t/\sqrt{3})e^{i(\omega_1t+\theta_1+\frac\pi2)},\\
\Omega_2(t)=-\sin(t/\sqrt{3})e^{i(\omega_2t+\theta_3-\theta_1+\frac\pi2))},
\ea\right.
\eqnl{OM2}
where $\omega_1=E_2-E_1,~\omega_2=E_3-E_2$. \\
 \\
Formula \r{OM2} coincides with formula \r{main-res} setting
$\varphi_1:=\theta_1+\pi/2$,  
$\varphi_2:=\theta_3-\theta_1-\pi/2$.\\  
We recall that in previous papers for the real-resonant problem, (see \cite{BCG}
for more details),
optimality was set as an assumption while here, it is obtained as a
consequence of the PMP. 

Moreover we get that the optimal cost and the probabilities $|c_j(t)|^2,$
$j=1,2,3$, are independent from
$\theta_1,\theta_3\in[-\pi,\pi]$.

\vspace{1cm} 
\noindent 
{\Large {\bf \underline{Appendix A: An interesting consequence of
the Cartan Decomposition}}}
\\
\bt 
Let $\l$ be a semi-simple Lie algebra and $\l=\k\oplus\p$ a Cartan
decomposition (i.e. it holds 
$[\k,\k]\subseteq\k,$ $[\p,\p]\subseteq\k,$ $[\k,\p]\subseteq\p$). Then
$[\k,\p]=\p$,  $[\p,\p]=\k$.  
\llabel{t-semisemi}
\et
\proof
We can  restrict ourself to the case in which $\l$ is simple, since a
semi-simple Lie algebra is a direct sum of simple ideals.\\\\
{\bf Claim 1:} $\p+[\p,\p]$ is an ideal in $\l$\\\\
Proof of Claim 1:
\bi 
\i  $[\p,\p+[\p,\p]]=[\p,\p]+[\p,[\p,\p]]$. Using the Cartan relations it
follows that the second term is contained in  $\p$. 
Hence
$[\p,\p+[\p,\p]]\subseteq \p+[\p,\p]$.  
\i   $[\k,\p+[\p,\p]]=[\k,\p]+[\k,[\p,\p]]$. Now, using the Cartan
relations we have:  $[\k,\p]\subseteq\p$ while
$[\k,[\p,\p]]=-[\p,[\k,\p]]-[\p,[\p,\k]]\subseteq [\p,\p]$, where we have
used the Jacobi identity.
Therefore
$[\k,\p+[\p,\p]]\subseteq \p+[\p,\p]$.  
\ei
{\bf Claim 2:} $\k+[\k,\p]$ is an ideal in $\l$\\\\
Proof of Claim 2:
\bi 
\i  $[\p,\k+[\k,\p]]=[\p,\k]+[\p,[\k,\p]]\subseteq[\k,\p]+\k$.
\i  $[\k,\k+[\k,\p]]=[\k,\k]+[\k,[\k,\p]]\subseteq\k+[\k,\p]$.
\ei
From the fact that $\l$ is simple (the only ideals are $0$ and $\l$) it
follows
$\p+[\p,\p]=\l\Rightarrow[\p,\p]=\k$, 
$\k+[\k,\p]=\l\Rightarrow[\k,\p]=\p$. \quadp\\\\
Notice that in general is false that $[\k,\k]=\k$. In particular it is
false for the Cartan decomposition of $SU(3)$ we used in this paper ($\k$
has a center). Theorem \ref{t-semisemi} is crucial to prove that,
in our case, abnormal
extremals are never optimal (see Appendix C).\\\\\\
\vspace{1cm} 
\noindent 
{\Large {\bf \underline{Appendix B: The $K+P$ Problem}}}
\\  
Here following \cite{jurd-book} and \cite{jurd-MCT}, we recall how to
write the PMP for a right-invariant
control system on a Lie group $G$.
We do it for a group of matrices only.

Let $\l$ and $\l^\ast$ be respectively the tangent and cotangent planes at
the identity of $G$. Consider the right-invariant system:
\bqn
\dot g=X(u) g,
\eqnl{dyn}
where $X(u)\in\l$ and $u$ belongs to the set of values of controls
$U$. To have a right-invariant optimal control problem, we must also assume a
right-invariant cost (that is a cost that does not depend on g, the coordinate
 on the group): 
\bqn \int_0^Tf(u(t))~dt,~~~\mbox{ where }f:U\to\R\mbox{ is
a
smooth
function.}  
\eqnl{cost}
Moreover we assume that the initial and final points belong to two given
smooth manifolds:
\bqn
g(0)\in M_{in},~~~~g(T)\in M_{fin}.
\eqnl{constr}
For each $\lam\in\R$ and $u\in U$ define the Hamiltonian :
$$
H_u(.):\l^\ast\to\R,~~~H_u(P):=P(X(u))+\lam f(u).
$$
We have $dH_u(P)\in(\l^\ast)^\ast$, ($P\in\l^\ast$) and since $\l$ has
finite dimension, we
can identify $(\l^\ast)^\ast\sim \l$ and consider $dH_u(P)\in \l$.
Let $B\in\l$ and denote with $ad^\ast_B(.)$ the operator from
$\l^\ast\to\l^\ast$ defined by:
$$
[ad^\ast_B(P)](A):=P\big(ad_B(A)\big)=P([B,A]),~~~\forall~A\in\l.
$$
The PMP (see \cite{pontlibro}) is a necessary condition for
optimality. For right-invariant control systems it reads (see
\cite{jurd-book, jurd-MCT}:
\bt
{\bf (PMP for right-invariant control systems)} Consider a couple
$(u(.),g(.)):[0,T]\to U\times G$ subjected to the dynamics \r{dyn} and to
the
constraints \r{constr}. If it minimizes the cost \r{cost}, 
then there exists a constant
$\lam\leq0$ and a never vanishing 
absolutely continuous function 
$P(.):t\in [0,T]\mapsto P(t)\in\l^\ast$ such that: 
\bqn
&&\frac{dg(t)}{dt}=dH_{u(t)}(P(t))g(t)\llabel{HS1}\\
&&\frac{dP(t)}{dt}=-ad^\ast_{dH_{u(t)}(P(t))}(P(t)).\llabel{HS2}
\eqn
Moreover: 
\bqn
&&H_u(t)(P(t))={\cal H}(P(t))~~
\mbox{  where  }~~~{\cal H}(P(t)):=\max_{v\in
U}H_v(P(t)).\llabel{HS3}\\
&&P(0)(T_{g(0)}M_{in}.g^{-1}(0))=0,~~~P(t)(T_{g(T)}M_{fin}.g^{-1}(T))=0
\mbox{~~ (transversality conditions) ~~}\llabel{transv-ap}
\eqn
\llabel{t-PMP}
\et
Here $g(t)\in G$ and $P(t)\in\l^\ast$ is the covector translated back to
the identity.
 $P(t)$ is related to the usual covector $p_{g(t)} \in T_{g(t)}^{\ast}G $ by:
  $P(t)=dR_{g(t)}^{\ast}p_{g(t)}$ or $p_{g(t)}=dR_{{g(t)}^{-1}}^{\ast}P_{g(t)}$

The couples $(u(.),g(.))$ satisfying conditions \r{HS1}, \r{HS2} and
\r{HS3}
with $\lam=0$ are called {\sl abnormal extremals}. Couples  $(u(.),g(.))$
corresponding to $\lam\neq0$ (in this case we can normalize  $\lam=-1/2$)
are called  {\sl normal extremals}. 

In the following we show that for normal extremals this 
Hamiltonian system
becomes completely
integrable if the problem is $\k\oplus\p$. Abnormal extremals may exist in
general, but they are never optimal as explained in Appendix C.

Consider  the right-invariant $\k\oplus\p$ problem defined in
Definitions
\ref{d-k+p}, \ref{d-k+p-problem}. Here $\l=\k\oplus\p$ can also be
noncompact
(in that case $\k$ must be the maximal compact subalgebra of $\l$, cfr.
Remark \ref{r-maximalcompact}). 

The Killing form defines a non-degenerate pseudo scalar product on $\l$.
This permits to identify $\l$ with $\l^\ast$ by:
\bqn
P\in\l^\ast\longleftrightarrow M\in\l\Longleftrightarrow
P(C)=Kil(M,C),~~\forall~ C\in\l.
\eqnl{idid}
Let us ``translate'' equation \r{HS2} for $M\in\l$.
For each $C\in\l$ we have:
\bqn
Kil\big(\frac{dM}{dt},C\big)&=&\frac{dP}{dt}(C)=
-[ad^\ast_{dH(P(t))}(P(t))](C)=-P(t)\left(  [dH(P(t)),C]  \right)\nn\\
&=&-Kil\big(M(t),[dH(P(t)),C]   \big)=+Kil\big([dH(P(t)),M],C \big),
\eqnn
where we have used the invariance of the Killing form under Lie Brackets: 
$Kil([A,B],C)=Kil(A,[B,C])$, that can be easily checked. The equation
for
$M(t)$ is then in the famous Lax-Poincar\'e form:
\bqn
\frac{dM(t)}{dt}=[dH(P(t)),M(t)],~~~M(t),~dH(P(t))\in\l.\llabel{poincare}
\eqn
Let $\{X_j\}$ be an orthonormal (right-invariant) frame for the $\p$ part 
of
$\l$, with respect to the
metric defined in Definition \ref{d-k+p-problem}. We have then 
$X(u)=\sum_j u_jX_j$ and $f(u)=u_1^2+...+u_{n_p}^2$ (here
$u=(u_1,u_2,.....u_{n_p})$ and $n_p$ 
is the dimension of the $\p$ subspace).
Moreover decompose $M=M_p+M_k$, where $M_p\in\p$ and  $M_k\in\k$.  
\bp
We have $dH(P(t))=X(u(t))=M_p(t)$.
\llabel{p-scema}
\ep
\proof The first equality can be obtained comparing equation $\dot
g(t)=X(u(t))g(t)$ with
equation \r{HS1} and using the fact that the differential of
the right translation is a linear isomorphism. Let us consider
the second one. 
From the maximum condition \r{HS3} with $\lam=-1/2$ and
$f(u)=u_1^2+...+u_{n_p}^2$,
we get: 
$$u_i(t)=
p_{g(t)}(t)(X_ig(t))=P(t)(X_i)=Kil(M(t),X_i)=Kil(M_p(t),X_i) 
\mbox{ because }
X_i=X_i(e)\in \p.$$
In the last equality we used the fact that $M_p$ and $M_k$ are
orthogonal
subspaces for the Killing form.
Therefore $M_p(t)=\sum_j Kil(M_p(t),X_j )X_j=\sum_j u_j(t)X_j=X(u(t))$.

\quadp\\\\
With Proposition \ref{p-scema}, the equation for $M$ become:
\bqn
\frac{dM_p}{dt}+
\frac{dM_k}{dt}=[M_p(t),M_p(t)+M_k(t)]=[M_p(t),M_k(t)].
\eqnl{Mpk}
Using the Cartan commutations relations \r{cartan-dec}, we have
$[M_p(t),M_k(t)]\subset \p$ and equation \r{Mpk} splits into:
\bqn 
\frac{dM_k}{dt}&=&0\Rightarrow M_k(t)=M_k(0),\nn\\
\frac{dM_p}{dt}&=&[M_p(t),M_k(0)]. \nn
\eqnn
Hence all the $\k$-components of the covector are constants of the motion.
Integrating
the equation for $M_p$, and setting $A_p:=M_p(0)$, $A_k:=M_k(0)$, we get:
\bqn
M_p( t)=e^{-A_k t}A_p e^{A_k t}=e^{-ad A_k t}A_p.
\eqnn
From equation \r{HS1}, with $dH(P(t))=M_p(t)$, we have:
\bqn
\frac{dg(t)}{dt}=\left(e^{-A_k t}A_p e^{A_k t}\right)g(t)
\eqnn
and the solution is equation \r{jurdjevic}:
\bqn
g(t)=e^{-A_k t}e^{(A_k+A_p)t} g(0)
\eqn 
Setting $A:=A_p+A_k$,  the transversality conditions
\r{transv-ap} reads in $\l$:

\bqn
Kil(A,T_{g(0)}M_{in}.g^{-1}(0))&=&0,\llabel{ttt-ap-1} \\
Kil(M(t),T_{g(t)}M_{fin}.g^{-1}(t))&=&0.\llabel{ttt-ap-2} 
\eqn

\brem {\bf (n-level case)} 
In the case in which the distribution is only a strict subspace of $\p$, 
as in the
$n$-level case ($n\geq4$), $dH(P(t))=X(u(t)$, but it is not equal to
$M_p(t)$ in general.
Then equation \r{Mpk} become:
\bqn
\frac{dM_p}{dt}+
\frac{dM_k}{dt}=[X(u(t)),M_p(t)+M_k(t)],
\eqnn
where now the right-hand side is not completely contained in $\p$. This
means that not all the 
$\k$-components of the covector are constants of the motion, and the proof
of the integrability of the Hamiltonian system given above fails.
\llabel{r-appendix}
\erem

\vspace{1cm} 
\noindent 
{\Large {\bf \underline{Appendix C: The Goh Condition}}}\\\\  
Let us consider any subriemannian metric over a manifold $M$, defined by
its orthonormal frame $\{X_i,~i=1,...,p\}$, completely nonintegrable. Then
a necessary condition for a strictly abnormal extremal (i.e.
an abnormal
extremal which is not normal at the same time) to be optimal is that it
satisfies the {\it Goh condition}.\\\\
{\bf Definition-Theorem} \textbf{(Goh-condition)} {\it Let $(p(t),x(t))$, $t\in[0,T]$ be the
Hamiltonian lift of the abnormal extremal $x(t)$ (see Theorem
\ref{t-PMP}). Then:
\bqn
(p(t)(X_i(x(t))\equiv0,~~p(t)\neq0,~~t\in[0,T]\llabel{goh-1},
\eqn
and a necessary condition for optimality of $x(.)$ is 
$\{p(t)(X_i),p(t)(X_j)\}(x(t))\equiv0,~~t\in[0,T],$ or:}
\bqn
p(t)\big([X_i,X_j](x(t))\big)\equiv0,~~t\in[0,T].\llabel{goh-3}
\eqn
This Theorem is a consequence of a (highly non trivial) generalized Maslov
index theory developped in \cite{AS}. For our right-invariant problem,
the relations \r{goh-1}, \r{goh-3} give:
$$p(0)(X_i(id))=0,~~p(0)\big([X_i,X_j])(id)\big)=0,$$ they imply with
Theorem \ref{t-semisemi} of Appendix A, that $p(0)(\l)=0$. Then $p(0)=0$.
This is a contradiction since $p(t)$ has to be nonzero, for all $t$.\\
\\
Hence, strictly abnormal trajectories are never optimal in our $\k \oplus 
\p$ problem.


\vspace{1cm} 
\noindent 
{\Large {\bf \underline{Appendix D: A technical computation}}}
\\
\bl $\mbox{Set } f_{a}=\cos(ta)\sin(t\sqrt{1+a^2})\frac{a}{\sqrt{1+a^2}}-\cos
(t\sqrt{1+a^2})\sin(ta), \mbox{ then } |f_{a}|\leq 1$.\\ Moreover, $|f_{a}|=1$ iff 
$\frac{|a|}{\sqrt{1+a^2}}=\left|\frac{1}{2k}+\frac{k'}{k} \right| <1$ ; $k\neq
0$
and $t=\frac{k\pi}{\sqrt{1+a^2}}$.\\
In particular, the smallest $|t|$ is obtained for $k=\pm 1$, 
$a=\pm \frac{1}{ \sqrt{3}}$, $t=\frac{\pm \pi \sqrt{3}}{2}$.
\el
\proof Set $\lambda=\frac{a}{\sqrt{1+a^2}}$, $\theta=t\sqrt{1+a^2}$, then
 \begin{eqnarray*}
f_{a}(t) & = & \lambda\cos(\lambda \theta) \sin(\theta) -
 \cos(\theta)\sin(\lambda \theta)\\
 & = & \langle (\lambda\cos(\lambda \theta),\sin(\lambda
 \theta)),(\sin(\theta),-\cos(\theta))\rangle \\
 & = & \langle v_{1},v_{2} \rangle .
 \end{eqnarray*}
  Both $v_{1}, v_{2}$ have norm $\leq $ 1 and $|f_{a}|\leq 1$.
  Hence, for $|f_{a}| = 1$, we must have $||v_{1}||=||v_{2}||=1$, $v_{1}=\pm v_{2}$.
   It follows that $\cos(\lambda\theta)=0$ and $\cos(\theta)
  =\pm 1$.
  Hence $\theta=k\pi$, $\lambda\theta=\frac{\pi}{2}+k'\pi$,
  $\lambda=\frac{1}{2k}+\frac{k'}{k}$. 
  Therefore, $\left|\frac{1}{2k}+\frac{k'}{k} \right|=\lambda<1$. Conversely,
  choose $k,k'$ meeting this condition, and $\theta=k\pi$. Then
  $\cos(\theta)=\pm 1$, $\sin(\lambda\theta)=\pm 1$, $f_{a}(t)=\pm 1$.\\
  Now, $|t|=\frac{k\pi}{\sqrt{1+a^2}}$, and the smallest $|t|$ is obtained for
  $k=\pm 1$ (if $k=0$, $\theta=0$ and $f_{a}(t)=0$). Moreover,
  $\left|\frac{1}{2k}+\frac{k'}{k} \right|<1$ is possible only for $(k,k')=(1,0)$ or
  $(1,-1)$ or $(-1,0)$ or $(-1,-1)$.\\
  In all cases, $|\lambda|=\frac{1}{2}$, $a=\pm\frac{1}{\sqrt{3}}$, and
  $t=\pm\frac{\pi\sqrt{3}}{2}$. \quadp



\fine